\documentclass{article}
\usepackage{maa-monthly}

\newtheorem{theorem}{Theorem}

\newtheorem*{definition}{Definition}
 
\newtheorem{proposition}{Proposition} 
\newtheorem{lemma}{Lemma} 
\newtheorem{fact}{Fact} 
\newtheorem{claim}{Claim}  
\newtheorem{conjecture}{Conjecture}

\newcommand{\proofClaim}{\noindent{\bf Proof (of the Claim). }}

\begin{document}

\title{Self Similarities of the  Tower of Hanoi Graphs  and  a proof of the Frame-Stewart Conjecture}
\markright{ Tower of Hanoi and Frame-Stewart conjecture}
%
{
\author{ 
Janez \v{Z}erovnik  \\
University of Ljubljana, FME,
A\v sker\v ceva 6 \\ 
and \\
IMFM, Jadranska 19,
SI-1000 Ljubljana, Slovenia\\
{\tt janez.zerovnik@fs.uni-lj.si;janez.zerovnik@imfm.si} 
}

\date{\today}

\maketitle

\begin{abstract} 
Considering the symmetries and self similarity properties of the corresponding labeled graphs,
it  is shown that the  minimal number of moves in the Tower of Hanoi game with  $p =4$ pegs and  $n \geq p$ disks
satisfies the  recursive formula 
$ F(p,n)   =  \min_{1\leq i \leq n-1} \{ 2F(p,i) + F(p-1,n-i) \} $
which proves the strong Frame-Stewart conjecture for the case $p=4$.
The method can be generalized to $p>4$.
\end{abstract}

\section{Introduction.}

The Tower of Hanoi game
with $p \geq 3$ pegs and  $n \geq 0$   disks is a popular game in which  disks have to be moved from one to another peg,
 obeying the rule that a larger disk can never be put onto a smaller one \cite{book}.
The case $p=4$ is sometimes called  Reve's puzzle. 
It is usual to study the game by regarding  the graph of possible positions  and legal moves among them.
The most prominent open problem is the Frame-Stewart conjecture
about  the  minimality  of  a  certain  algorithm  for  the  distance  between  
so called perfect    states.
The known solution  to the Problem 3918 \cite{problem}  
that has been reinvented several times (see \cite{book} for historical facts)
appears to be a very natural one, 
but the proof of optimality  in general case  is not known until today.
The two solutions proposed by Frame and Stewart  \cite{FrameStewart1,FrameStewart2}
 are known to be equivalent \cite{Klavzar}
and the number of steps needed for $p$ pegs and  $n\geq p$ disks 
is given by the   recursive formula
\begin{equation}
F(p,n)   =  \min_{1\leq i \leq n-1} \{ 2F(p,i) + F(p-1,n-i) \} .  \label{FrameStewart}
\end{equation}
As it is trivial that $F(p,n) =2n-1$ for $n<p$, 
and it is well-known that $F(3,n) = 2^n -1$, 
  formula (\ref{FrameStewart})   determines $F(p,n)$ for all $p\geq 3$ and $n\geq 0$. 
The literature on the topic is enormous, for example   there are 352 references in \cite{book}.
Very recently,  a   proof of optimality  for the case $p=4$ appeared in \cite{Bousch}.

Here we first observe some self similarity properties of the  labeled graphs that are 
isomorphic to Hanoi graphs. 
The symmetries can rather  naturally be observed from drawings that, to our surprise, 
are used very rarely. In fact, in an attempt to Google search 
for drawings of Hanoi graphs,  the present author found a similar drawing in only one paper \cite{Zhang}, which seems to be unpublished students' homework(!).
The structure of the Hanoi graphs allows to prove a lemma about existence of shortest paths of certain structure that in turn provides a proof  that the Frame-Stewart algorithm is optimal. This solves  the famous Frame-Stewart conjecture.
 
The rules of the game are simple.
There are $p$ pegs and $n$ disks, that  all differ in size (diameter).
In the beginning, all the disks are at one peg, (legally) ordered by size.
It is not allowed to put a larger disk onto a smaller one at any time.
In one move,  a disk that is on the top at one peg is  put to  the top  at another peg.
The task is to reach a state in which  all the disks from original peg  are at  another peg so that the number of moves is minimal.

 Only basic notions of graph theory will be used. 
A graph is a pair of sets  $G = (V(G),E(G))$, where $V(G)$ is an arbitrary set of vertices, 
and $E(G)$ is a set of pairs of vertices. 
Usual notation is $e=uv$ meaning that edge $e$ connects vertices $u$ and $v$.
We work with labeled graphs, i.e. graphs   with a labeling function that assigns a label to each vertex.
Two (unlabeled) graphs  $G$ and $H$ are isomorphic when  there is a bijection (called isomorphism)
$\alpha : V(G) \to V(H)$ such that $u$ and $v$ are connected in $G$ if and only if 
$\alpha(u)$ and   $\alpha(v)$ are connected in $H$.
A walk  is a sequence of vertices and  edges $v_0 e_1 v_1 e_2 v_2 ...e_k v_k$ 
such that  $e_i = v_{i-1}v_i$. 
A walk is also determined either by the  sequence of vertices  $v_0  v_1 v_2 ...v_k$ 
or by the sequence of edges  $ e_1  e_2  ...e_k $.
 The length of a walk is the number of edges on it.
A path is a walk in which all vertices are distinct.     
The distance between two vertices is the length of a shortest path.
A subgraph $H$ of a graph $G$ is isometric subgraph, if the distance between any two vertices in $H$ is equal to the distance  in $G$.
For notions not recalled here see, for example \cite{knjigaTG}.

The rest of the paper is organized as follows. 
In the next section, we construct labeled graphs $G^{(p)}_n$  corresponding to the Tower of  Hanoi   game 
with $p$ pegs and $n$ disks.
In Section 3, a lemma about the self similarity of these graphs is proved and some related facts are given.
Sections 4  and 5  provide proof of  existence of shortest paths of certain form 
which in turn gives a lower bound on length of the shortest paths between certain vertices.
This implies  the main result, a proof of the  Frame-Stewart conjecture for four pegs  that appears in Section 6.

%
\section{The construction of labeled graphs $G^{(p)}_n$.}\label{casekpegs}
%
 
Before giving the general definition, we start with the $p=4$  pegs example.
Labels of  the graphs $G^{(4)}_n$  will be words of length $n$ over the 
four letter  alphabet  $ {\cal A}_4= \{A,B,C, D\}$.

Let  $G^{(4)}_1$ be the tetrahedron graph (complete graph on 4 vertices), 
and the vertices labeled with $A,B,C,$ and $D$. 

Given $G^{(4)}_n$, we construct $G^{(4)}_{n+1}$ as follows. 
Take four copies of  $G^{(4)}_n$:
in the first copy, denoted by $AG^{(4)}_n$ replace each label $*$ with label $A*$
(i.e. add an $A$ at the beginning of each label).
Similarly,  in $BG^{(4)}_n$ replace each label $*$ with label $B*$,   
 in $CG^{(4)}_n$ replace each label $*$ with label $C*$,   and 
 in $DG^{(4)}_n$ replace each label $*$ with label $D*$. 

Finally connect some  pairs of vertices  from different copies of $G^{(4)}_n$ with edges 
by the following rule.
Let $X,Y \in {\cal A}_4$, $X\not=Y$.
Two vertices $u\in V(XG^{(4)}_n)$ and $v\in V(YG^{(4)}_n)$
are connected if and only if the labels of $u$ and $v$  
without the first letter are equal and   are words over ${\cal A}_4 - \{X,Y\}$.
(i.e. the labels do not contain any $X$ or $Y$).
The graphs $G_1^{(4)}$ and  $G_2^{(4)}$  are drawn on Fig.  \ref{G1inG2} and the graph  $G_3^{(4)}$  is on Fig.  \ref{G3classes}. 
 
In other words,  from  the definition it follows that two vertices from different copies of $G^{(4)}_n$
are connected  in $G^{(4)}_{n+1}$
  exactly when their labels  are equal  (in the last $n$ letters)
and the labels are words over an alphabet of two letters.

\begin{figure}[ht]
 \centering
  \leavevmode
  \includegraphics[height=7cm]{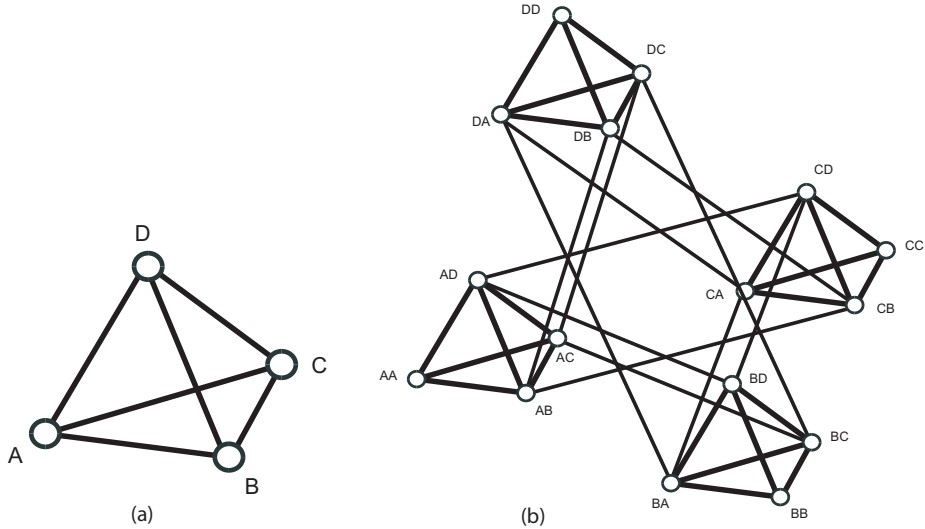} 
 \caption{The graphs $G_1^{(4)}$ and  $G_2^{(4)}$ .} 
 \label{G1inG2}
\end{figure}

\noindent  {\bf Remark. }
The graph  $G^{(4)}_n$ is isomorphic to the Hanoi graph of the game with 4 pegs and  $n$ disks.
Just interpret the labels naturally as: $i$-th letter in the label is the position of the $i$-th disk.
For example, the first letter gives the position of the largest disk.

We now turn to general definition, for arbitrary $p \geq 3$. 
A generalization to  $G^{(p)}_n$ is straightforward : 

\begin{definition}
 Let $G^{(p)}_1$ be a complete graph with $p$ distinct labels on vertices, say using letters from alphabet ${\cal A}_p$.
Construct  $G^{(p)}_{n+1}$   using $p$ copies of  $G^{(p)}_n$ as above, i.e.
in each copy of $G^{(p)}_n$    use a different letter as a prefix for labels.
As before, connect two vertices $u\in V(XG^{(p)}_n)$ and $v\in V(YG^{(p)}_n)$
  if and only if the labels of $u$ and $v$  
without the first letter are equal and   are words over ${\cal A}_p- \{X,Y\}$
(i.e. the labels do not contain any $X$ or $Y$).
\end{definition}

\noindent  {\bf Remark. }
By construction, the graph  $G^{(p)}_n$  is  the graph of the Tower of  Hanoi game with  $k$ pegs and $n$ disks.
(As unlabeled graph, it is therefore  isomorphic to the Hanoi graph,  $H_n^p$    in   notation of  \cite{book}.) 
The Frame-Stewart conjecture  says that the minimal number of moves  between two perfect states is determined by recursion  (\ref{FrameStewart}) 
which    is equivalent to statement that  this is the  length of a shortest path  in  $G^{(p)}_n$  between two perfect vertices,
for example  vertices  with labels $A^n = AA...A$ and $B^n = BB...B$.

%
\section{Self Similarity.}\label{SelfSimilarity}
%

\begin{figure}[htb]
 \centering
  \leavevmode
  \includegraphics[height=8cm]{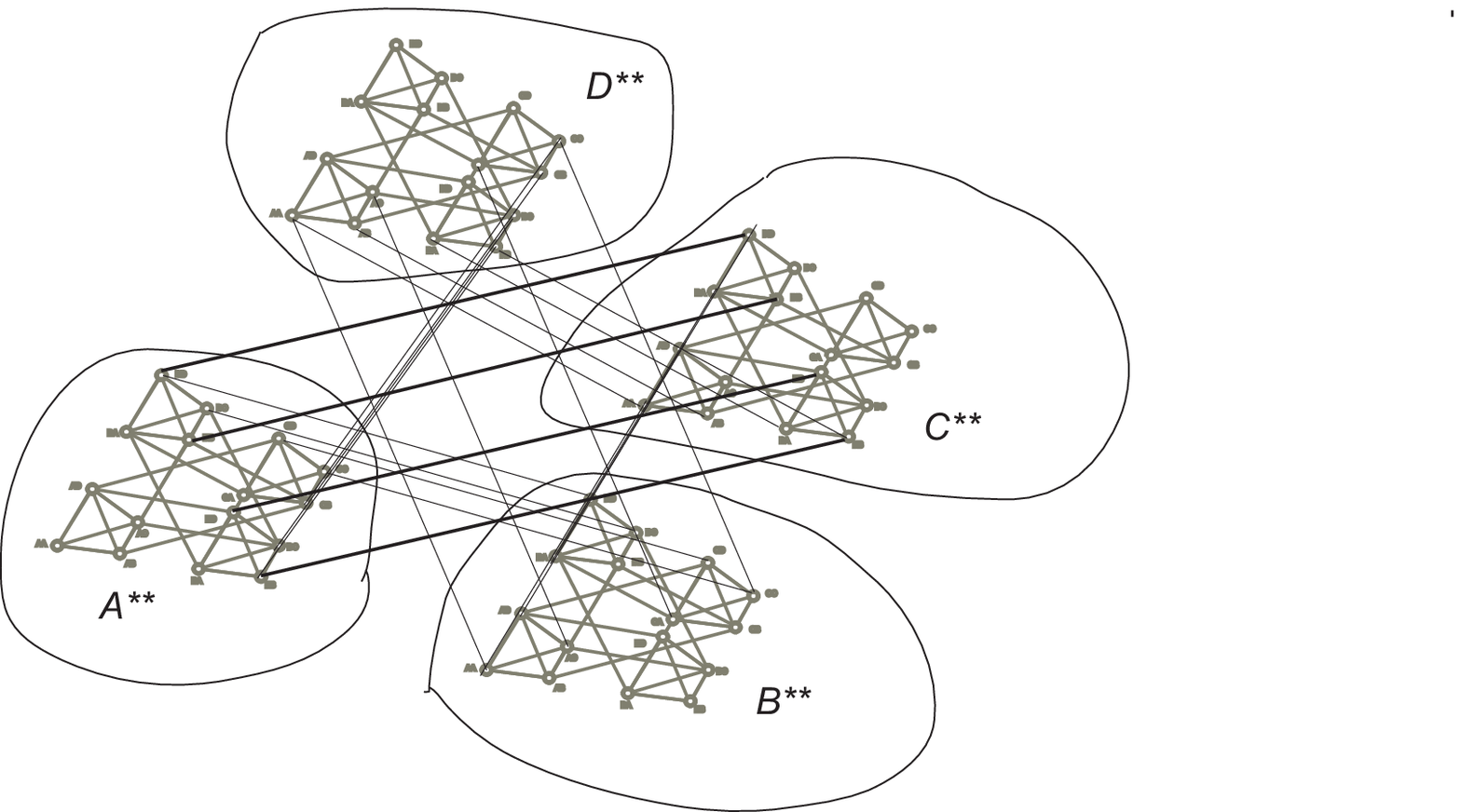} 
 \caption{The graph $G_3^{(4)}$; with the four equivalence classes ($i=1$).}  
 \label{G3classes}
\end{figure}

It is well known that the Hanoi graphs are highly symmetric.
Here we first  emphasize    one of the properties that motivated  the  idea used in the argument given  below.
Some more properties are listed below for later reference.  
We do not give detailed proofs as we believe that these results are not new, and are recalled here  for completeness of presentation.

Given  $G^{(p)}_n$ and $1\leq i\leq  n$  define the  equivalence relation $R_i$ 
on the set of vertices  $V(G^{(p)}_n)$ 
 as follows:
two vertices are equivalent when their labels coincide in the first $n-i$ letters.
Define the  graph   $G^{(p)}_n  / R_i$  on  equivalence classes (as vertices) by connecting two 
equivalence classes if there is an edge   $G^{(p)}_n $  that connects a pair of vertices from the two classes.
The definition directly implies:

\begin{lemma}
$G^{(p)}_n  /  R_i $ is isomorphic to $G^{(p)}_{n-i}$.
\end{lemma}

Having in mind this structure, we will (for fixed $i$)
make a distinction between the  edges within equivalence classes   
and the edges connecting different classes. 
The later will be called {\em bridges} and the edges within equivalence classes will be referred to as {\em local edges}.
Thus bridges correspond to moves of the largest $n-i$ disks and 
local edges to moves of the smallest $i$ disks.
(See Fig.  \ref{G3classes}  and   Fig.  \ref{G3edges}, where the bridges  between classes A** and C** in $G_3^{(4)}$ are shown.)

\begin{figure}[htb]
 \centering
  \leavevmode
  \includegraphics[height=8cm]{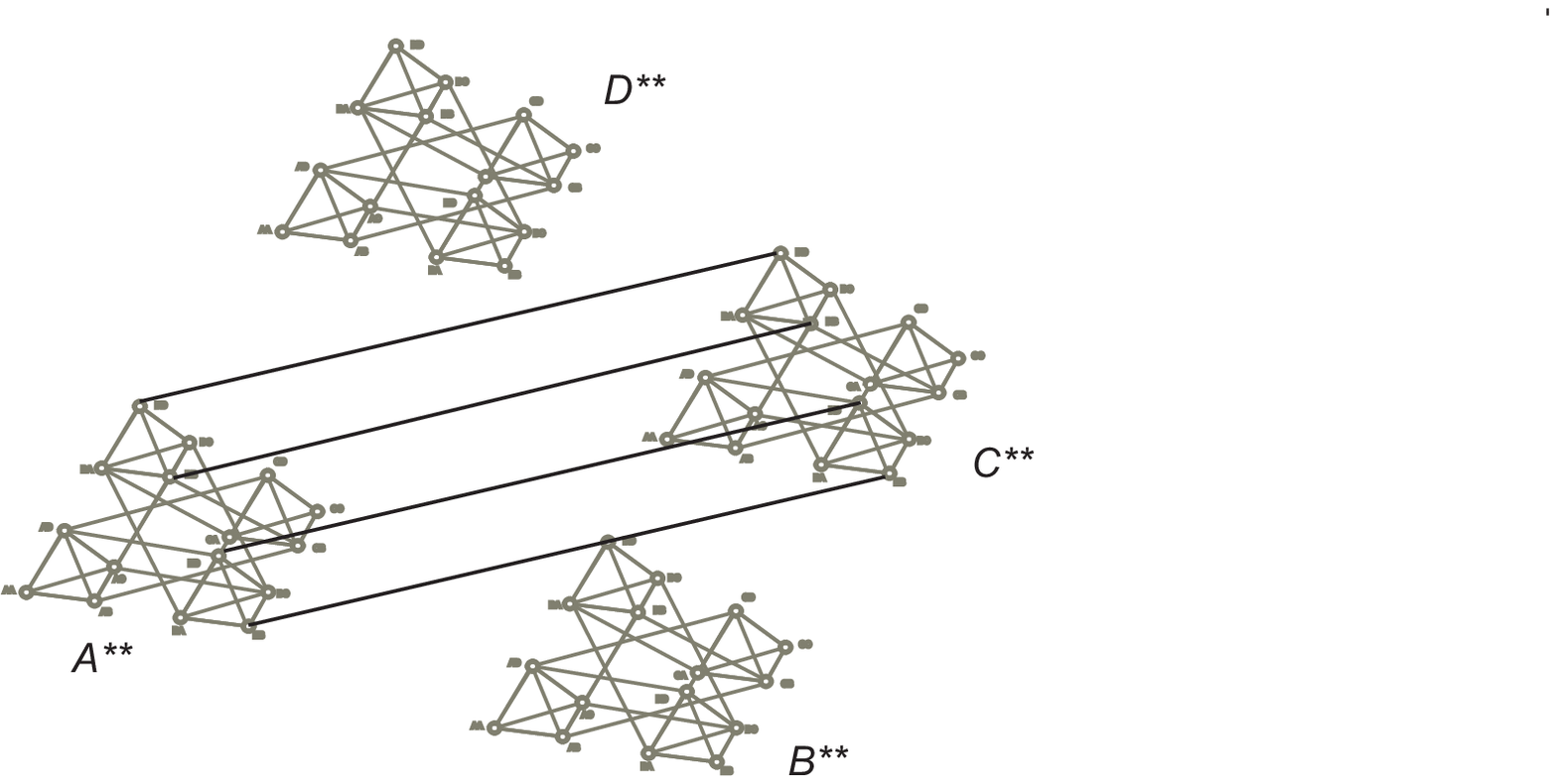}
 \caption{The four equivalence classes ($i=1$) of the graph $G_3^{(4)}$ and all edges (bridges)  between two classes.}  
 \label{G3edges}
\end{figure}

By definition of $R_i$, each equivalence class of $R_i$ induces a subgraph of $G^{(p)}_n$ 
that is isomorphic to $G^{(p)}_{i}$.
Furthermore, recall that any bridge connects two vertices (from different equivalence classes) 
with the same labels (i.e. with labels that match in the last $i$ letters).

%

We state two more properties  for a later reference. 
The proofs follow directly from the definitions. 
Maybe even more natural argument is to consider the meaning in terms of the game, 
namely (1) fixing positions of the largest $n-i$ disks clearly results in a game with $i$ disks and 
(2) putting the smallest $i$ disks to one of the pegs forbids the moves to that peg, hence 
exactly $p-1$ pegs are free to move the largest  $n-i$ disks at.

\begin{lemma} \label{LemaIsometry}
Let $W$ be a arbitrary word of $n-i$ letters over alphabet ${\cal A}_p$. Then 
 $WG^{(p)}_{i}$  are isometric subgraphs of  $G^{(p)}_n$.
\end{lemma}

\begin{lemma}  \label{Lema33}
The vertices 
 $WG^{(p)}_{i}$ of   $G^{(p)}_n  /  R_i $  
where $W$  is an  arbitrary word of $n-i$ letters over alphabet ${\cal A}_p - X$
induce a subgraph that    is isomorphic to  $G^{(p-1)}_{n-i}$.
\end{lemma}

We conclude the section with a couple of facts that will be useful later. 
As above, the proofs are not difficult, for example by considering the meaning of the distances in terms of the game.
In short, the first claim follows by obvious symmetry (just replace the role of $A$ and $B$). 
The second claim is obvious because if  we can solve the bigger task, then we also can solve the easier task (and need not move the largest disk). 
Details are left to the reader.

\begin{fact} \label{Fact1} 
Let $A$ and $B$ be arbitrary letters from alphabet   ${\cal A}_p $. 
Let $W$ be a word of $i$ letters over alphabet   ${\cal A}_p- \{A,B\}$.
Then the distance between vertices with labels $A^i$  and $W$ is  equal to  the distance between vertices with labels $B^i$  and $W$.
Hence, there is no shortest path connecting vertices with labels $A^i$ and $W$ that meets $B^i$  in $G_i^{p}$.
\end{fact}

\begin{fact} \label{Fact2}  
Let $W$ be a word of $i$ letters over alphabet   ${\cal A}_p- \{A\}$.
Then there exists   $C\not= A $
such that   the distance between vertices with labels $C^i$ and $W$ is strictly smaller  than 
 the distance between vertices with labels $A^i$ and $W$.
\end{fact}

\noindent  {\bf Remark. }
Idea of proof of Fact \ref{Fact2} :
On any path  $P$   from $W$ (i.e. vertex   with  label $W$) to $A^i$, 
the largest disk will be moved to peg $A$, so the path $P$ can be written as $W \to W^{*} \to AW^{**} \to A^i$, 
where $W^{*}$ and $ AW^{**} $  are two labels that differ only in the first letter. Let the first letter of  $W^{*}$  be $C$.
Then there is a path   $P^{*}$  from $W$ to $C^i$,    i.e. $W \to W^{*}  = CW^{**} \to C^i$,  and 
 $P^{*}$ is shorter than $P$.

%
\section{Shortest paths.}\label{ShortestPaths}
%
 
Let $n\geq 2$, as the case $n=1$ is   trivial.
Let ${\cal A}_p$ be an alphabet of $p$ letters, and $A,B \in  {\cal A}_p$.

Let $P$ be a shortest   path  from vertex $a$ with label $A^n$  to vertex $b$ with label  $B^n$.
Below, in  Proposition   \ref{LemaExistence}  we will assume 
 that on the path $P$ there is a vertex that has a label of the form 
$A^{n-i}X^i$, where $X \in {\cal A}_p -\{A,B\}$  (i.e. $X$ is not $A$ nor $B$).
We call such a vertex {\em special}. 
It may seem obvious that there is a special vertex on every shortest path from  $a$ and $b$. 
However, this is not the case as pointed out by Ciril Petr and Sandi Klav\v zar.  
A slightly weaker,  but still sufficient,   statement can be proved

\begin{lemma}  \label{lemaspecial}
Let $p=4$. 
There is  a shortest path   with at least one  special vertex on every path between vertices   $a$ and $b$  with labels $A^n$   and  $B^n$.
\end{lemma}

The proof  of Lemma   \label{lemaspecial} is postponed to the next section.  It is based on induction, in which the small cases are, due to large number of them, 
rather tedious task.  
Alternatively, it can be checked by computer using a straightforward  application of a shortest path algorithm.
We believe that the same technique can be used for $p>5$ and conjecture

\begin{conjecture}   \label{special}
There is  a shortest path   with at least one  special vertex on every path between vertices   $a$ and $b$  with labels $A^n$   and  $B^n$.
\end{conjecture}

Now  we will show that there is a shortest path with certain structure. 
More precisely,  

\begin{proposition} \label{LemaExistence}
Let $P$ be a shortest  path connecting vertices  $a$ and $b$ with labels  $A^n$   and $B^n$  and let $1\leq i\leq n$.
If there is a special vertex on $P$ with label of the form  $A^{n-i}X^i$,  where $X \in {\cal A}_p -\{A,B\}$ 
then   there is a shortest  path $Q$ from  $a$  to $b$  in $G^{(p)}_n$, 
which is a concatenation of three subpaths  $Q_1$,  $Q_2$,  $Q_3$
such that $Q_1$ and $Q_3$ only use local edges and $Q_2$ only uses bridges. 
\end{proposition}

\proof
Consider a shortest   path $P$  from  $a$ to $b$ (vertices with labels $A^n$   and $B^n$).
Let $i  = i(P)$ be maximal with property that there is a special vertex $s = s(i)$ with label  $A^{n-i}X^i$ on the path $P$.

Denote by $P_1$ the first part of $P$, from  $a$ to $s$, 
and observe that because of Lemma \ref{LemaIsometry} we can, without loss of generality, assume that there are only local edges on $P_1$. 
As $P$ is a shortest path, $Q_1 =P_1$ must be a shortest  path from  $a$ to $s$ within the subgraph $A^{n-i}G^{(p)}_i$.
 
On path $P$ there must be at least one bridge after $P$ meets $s$
because labels 
 $A^n$ and $B^n$ differ in the
 first  $n-i \geq 1$ letters that can only be changed when traversing bridges.
First we show that there is a shortest path  such that the edge used after visiting $s$ is a bridge:

\begin{claim} \label{Claim1} 
Let $P$ be a shortest path  from   $a$ to $b$  (with labels  $A^n$  and $B^n$) that meets vertex $s = s(i)$ with label  $A^{n-i}X^i$, $i  = i(P)$.
Then there is a   path $P^\prime$ of the same length such that the first edge after visiting $s$ is a bridge.
\end{claim}

\proofClaim  
If $P$ has the property claimed then $P^\prime  =P$ and we are done.
Now assume that $P=P_1P_2 f P_3$ where $P_1$ is a shortest path from  $a$  to $s$,
$P_2$ is a subpath of local edges from $s$ to $w$, $f$ is a bridge, and $P_3$ is the rest of $P$.

We distinguish two cases.
First, let $w$ be a vertex with label $A^{n-i}W$ where $W$ is a word using letter $X$.
This implies that $f$ is a bridge that moves a disk from peg $A$ to a peg that is not peg $X$, 
say $f$ moves $n-i$-th disk from $A$ to $Y\not=X$.
This implies that  we can replace subpath $P_2 f$ with $f^\prime P_2^\prime$, 
where  $f^\prime$ is the edge connecting vertex $s$ (with label $A^{n-i}X^i$) and  vertex $s^\prime$ that has label  $A^{n-i-1} Y X^i$
and  $P_2^\prime$ is a copy of  $P_2$  in the subgraph   $A^{n-i-1} Y G^{(p)}_i$.
(Formally, the path    $P_2^\prime$ is constructed from   $P_2$ by replacing  every vertex on $P_2$ (with label   $A^{n-i}*$)
with the vertex with label $A^{n-i-1}Y*$.)

Now assume that  $w$ is  a vertex with label $A^{n-i}W$ where $W$ is a word without  letter $X$.
In this case  the subpath $P_1 P_2$  (within   $A^{n-i}   G^{(p)}_i$)
starts at $a$ (with  label $A^n$),  visits first vertex with label  $A^{n-i}X^i$ and then reaches  $w$. 
However, there is a strictly shorter path  from  $a$ to $w$  
which contradicts the assumption that $P$ is a shortest path  (recall Fact \ref{Fact1}).
\qed (Claim)

Now we will prove that there is no need to use local edges between bridges, i.e. that there is a shortest path which, after visiting $s$,
first traverses all bridges, and then traverses only local edges. 
The next claim  shows that we can always  conveniently  change position of the next displaced bridge. More precisely, 

\begin{claim} \label{Claim2}   
Let $P$ be a shortest path of the form $P= P_1 P_2 L f P_3$, 
where $P_1$ is a shortest path  from $a$   to $s$,
$P_2$ is a path of bridges from $s$ to $ t$, 
$L$ is a local path, $f$ is a bridge, and $P_3$ is the rest of $P$.
Then there is  shortest path $P^\prime = P_1 P_2 f^\prime L^\prime P_3$, 
where  $f^\prime$ is a bridge and $L^\prime$ is a local path of the same length as $L$.
\end{claim}

\proofClaim  
The local path $L$ is a path from vertex  $t$ to a vertex, say $w$.
Let the label of $t$ be   $ZX^i$ where $Z$ is a word of length $n-i$ that does  not contain any $X$.
Then the label of $w$ must be of the form $ZW$.

We distinguish two cases.
First, assume the word $W$ contains at least one letter $X$. Therefore bridge $f$ that moves one of the large disks may not move  a   disk to peg $X$.
(More precisely, 
edge $f$ connects $w$ and $u$, and the label of $u$   must be  $Z^\prime W$, where the words   $Z^\prime$ and $Z$ differ in one letter at  one position, and 
neither   $Z^\prime$ nor  $Z$  contain any letter $X$.)

Hence we can define $f^\prime$ to be the bridge that connects $t$ (with label $ZX^i$) and the vertex with label  $Z^\prime X^i$, denote it $r$.
Furthermore, let $L^\prime$ be a copy of $L$ in $Z^\prime  G^{(p)}_i$ that connects vertices $r$ and $u$. By construction, paths  $f^\prime L^\prime$ and $Lf$
both connect $t$ and $u$ and are of the same length, as needed.

Now assume that $w$ is a vertex with label $ZW$ where $W$ is a word without letter $X$. In this case, we can construct a shorter path from $a$ to $w$, which contradicts the assumption that $P$ is a shortest path. The argument is as follows.

First, observe that, because there is no $X$ in $W$, the first letter of $W$ is $Y\not= X$.
(In other words,  the $i$-th smallest disk is on peg $Y$ that is not peg $X$.)

Second, by symmetry, existence of the path $P_2$  (connecting vertices with labels  $A^{n-i}X^i$  and $ZX^i$) implies the existence of a path $P_2^\prime$
of the same length that connects vertices with labels   $A^{n-i}Y^i$  and $ZY^i$.
(Just replace any occurence  of $Y$  in $P_2$ with $X$  (and vice versa) to get $P_2^\prime$.)

Furthermore, by symmetry, there is a path $P_1^\prime$  (of the same length as $P_1$ that connects vertex $a$ and the vertex with label  $A^{n-i}Y^i$.

Finally, recalling Fact \ref{Fact2}, in $ZG_i^{p}$  there is a path $L^\prime$ from the vertex with label $ZYY^{i-1}$ to vertex $w$ that is strictly shorter than $L$.

Summarizing, the path  $P_1^\prime P_2^\prime L^\prime$ from $a$ to $w$ is shorter than $ P_1 P_2 L$, contradicting the minimality of $P$.
\qed (Claim)
  
By inductive application of the last Claim we prove that any shortest path $P$ can be replaced by a shortest path $Q=  Q_1 Q_2 Q_3$
where $Q_2$ is a path of bridges while $Q_1$ and $Q_3$ are path of local edges. 
\qed

%
\section{Special vertices - partial proof of Conjecture \ref{special}.}

In this section we will prove  Lemma  \ref{lemaspecial}. 
The proof is by induction. 
We emphasize that  the inductive step is  proved  for general $p$ 
while  we are only able to prove the base step(s) for each $p$ separately.

We begin by introducing some more notation. 
In the graph   $G_{n+1}^{(p)}$,
we consider shortest paths from $a$ to the set of vertices on the {\em border} of  $AG_{n}^{(p)}$, i.e. to vertices  
from which there are edges going out to other subgraphs  $XG_{n}^{(p)}$, $X\not=A$.
The border  $S_n$  is, by definition  of   $G_{n+1}^{(p)}$,  a union of 
$S_n^{AX}$,   where   $S_n^{AX} $ is a set of vertices with labels that do not use letters $A$ and $X$ :
$ \displaystyle{ S_n  =  \bigcup_{X\in {\cal A} - A}   S_n^{AX} } $.
As $AG_{n}^{(p)}$ is just a copy of $G_{n}^{(p)}$ we can recursively define $S_{n-1}$ in  $G_{n}^{(p)}$,   
and because $AG_{n}^{(p)}$ is a subgraph of  $G_{n+1}^{(p)}$, this also gives a definition of $S_{n-1}$ in  $G_{n+1}^{(p)}$.  
Thus we define the sets  $S_i$  in  $G_{n+1}^{(p)}$    for  $i=1,2,\dots, n$ (see Figure   \ref{Boundaries}). Obviously, 

\begin{fact} \label{InductionStep}
 Any shortest path from $a$ to $v \in S_{n}$   meets $S_{n-1}$.
\end{fact}

\begin{figure}[htb]
 \centering
  \leavevmode
  \includegraphics[height=8cm]{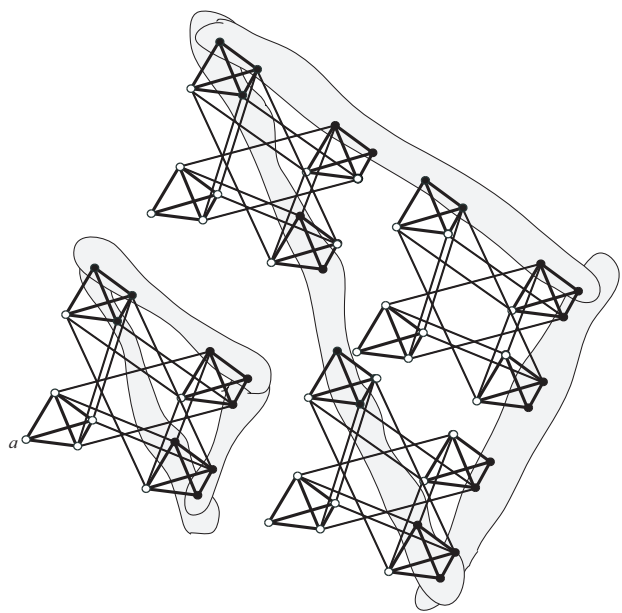}
 \caption{  Subgraph $A^{n-3}G_3^{(4)}$    with  emphasized  vertex sets $ S_2$   in $S_3$.  Edges corresponding to the moves of the 3rd disk are not drawn.}  
 \label{Boundaries}
\end{figure}

As already mentioned, there are examples, in which there  is a unique shortest path from $a$ to some vertices on the border. 
We conjecture that this is only possible for small $n$.

\begin{conjecture}   \label{SmallConjecture}
For any $p \geq 3$ there is  $n_0 \geq 2$ such that  for any $v \in S_{n_0}$  there is 
a  shortest path from $a$ to $v$   that   meets  a special vertex. 
\end{conjecture} 

For $p=4$, the validity of  conjecture can be   proved by regarding  $S_6$ in  the graph $G_7^{(4)}$.
In fact the  present  author performed    a  tedious case analysis "by hand" using a tool for 
drawing shortest paths in  Hanoi graphs by Igor Pesek.
As a curiosity, let us mention that the only vertices on the boundary $S_5$ that have a unique shortest  path to $a$
are the six vertices  with labels  $A...ACCDCD$, $ A...ADDCDC$, (and, by symmetry,  $A...ACCBCB$, $A...ABBCBC$,  $A...ABBDBD$, $A...ADDBDB$), 
two of them  passing the vertex with label $A...AAAAAB$. 
On $S_6$, there is no vertex with a unique  shortest path to vertex $a$, and in particular, in each case at least one of the paths   avoids vertex   with label $A...AAAAAB$. 

We thus know that

\begin{fact} 
For any $v \in S_6$  in  the graph $G_7^{(4)}$, there is  a shortest path from $a$ to $v$   that   meets  a special vertex.
\end{fact} 

Hence, 
by induction step (Fact  \ref{InductionStep}),  for any   $v \in S_{n}$ in  the graph $G_{n+1}^{(4)} $ 
there is a  shortest path from $a$ to $v$   that meets  a special vertex, for any $n\geq 7$.

Recall that  any shortest path  between $a$ and $b$  (with labels  $A^n$   and  $B^n$)  in  $G_{n }^{(p)}$
meets $S_{n-1}$, more precisely   $S_{n-1}^{AB}$.
Provided   validity of  Conjecture \ref{SmallConjecture}   (that is proved above for case $p=4$),
we have the existence of a shortest path  between  $a$ and $b$  (with labels   $A^n$   and  $B^n$)   
that meets a special vertex. 
Furthermore, for $p=4$ and $n<7$ it is straightforward to construct  shortest  paths that meet a special vertex. 
This concludes the proof of   Lemma  \ref{lemaspecial}.

\medskip
\noindent{\bf Remark.}
Clearly, the cases $p>4$ can be handled along the same lines. The size of graphs however is probably too large to  check without computer assistance.
 
%
\section{The main result.}

  Proposition   \ref{LemaExistence}  implies that there is a shortest path  between  vertices  $a$  and  $b$  (with labels  $A^n$   and   $B^n$)  
in which bridges are all sandwiched together between two local paths.
Hence, it suffices to consider the paths of this form.

\begin{lemma}\label{LemaDistances}  
Assume validity of Conjecture \ref{special}.
Let $ P$ be a path connecting vertices  with labels $A^n$   and $B^n$,
and    with  a special vertex on $P$  with label of the form  $A^{n-i}X^i$,  where $X \in {\cal A} -\{A,B\}$.
  Assume   $ P$   is a concatenation of three subpaths  $ P_1$,  $ P_2$,  and $ P_3$
such that $P_1$ and $P_3$ only use local edges and $P_2$ only uses bridges. 
Then the length of $ P$ is     $F(p,i) + F(p-1,n-i) + F(p,i)$.
\end{lemma}

\proof
Let the special vertex $s$ on $P$ be of the form   $A^{n-i}X^i$.
Then the last $i$ letters  of labels of  all the vertices on $P_2$ are $X^i$, 
which implies that $P_2$ never meets the copies $WG_{i}^{(p)}$
where $X$ appears in the word $W$.
Recall that  by Lemma \ref{Lema33},  the subgraph  induced on vertices  $WG_{n-i}^{(p)}$
where $W$ is a word of length $n-i$ over alphabet  ${\cal A}_p - X$  
is isomorphic to $G_{n-i}^{(p-1)}$.

The lengths  of   subpaths  $  P_1$,  $  P_2$,  $  P_3$  are thus bounded by  
$F(p,i)$,  $ F(p-1,n-i)$,   and  $ F(p,i)$, respectively. 
Hence the statement of the Lemma follows.
\qed 

\begin{theorem}  \label{main} 
Assume validity of Conjecture \ref{special}.
The distance between   vertices with labels $A^n$   and   $B^n$  in   $G^{(p)}_n$   is 
$  \min_{1\leq i \leq n} \{ 2F(p,i) + F(p-1,n-i) \}$.
\end{theorem}
 
\proof
Consider a shortest path $P$  between vertices with labels  $A^n$  and  $B^n$.
There must be a special vertex on $P$, 
and, by Lemma \ref{LemaExistence}, a shortest path $Q$  
for which the length is given by   Lemma \ref{LemaDistances}.
Recall  \cite{Klavzar,book} that  there are algorithms for which the number of moves is given by Eq. (\ref{FrameStewart}).
Therefore  the length  of  a  shortest path $P$ 
is of the from   $2 F(p,i(P)) + F(p-1,n-i(P))$, as claimed.
\qed 

Theorem \ref{main}    and Lemma \ref{lemaspecial}   impliy 

\begin{theorem}   \label{FS}
The strong Frame-Stewart conjecture is true for $p=4$.
\end{theorem}

Recall that for any $p$,  validity of    Conjecture \ref{special} implies the strong  Frame-Stewart conjecture for that $p$.

\begin{acknowledgment}{Acknowledgment.}
The author  wishes  to thank Igor Pesek for providing a tool for drawing shortest paths in Hanoi graphs, and  to Sandi Klav\v{z}ar and Ciril Petr
for useful comments  on an earlier version of this paper. The work was in part supported by ARRS.
\end{acknowledgment}

\end{document}